\documentclass[a4paper,11pt,twoside,english]{article}
\usepackage{babel}
\usepackage{latexsym,amsfonts}
\usepackage{amsthm}
\usepackage{amsmath}
\usepackage{mathtools}
\usepackage{paralist}

\usepackage{lastpage}
\usepackage{fancyhdr}
\fancypagestyle{plain}{

  \fancyhead[L,C]{}
  \fancyhead[R]{\today}
  \fancyfoot[R,L]{}
  \fancyfoot[C]{\thepage \ of \pageref{LastPage}}
}

\title{Summands in locally almost square and locally octahedral spaces}
\author{Jan-David Hardtke}
\date{}

\providecommand{\sm}{\setminus}
\providecommand{\ssq}{\subseteq}

\providecommand{\N}{\ensuremath{\mathbb{N}}}

\providecommand{\R}{\ensuremath{\mathbb{R}}}

\providecommand{\U}{\ensuremath{\mathcal{U}}}
\providecommand{\eps}{\ensuremath{\varepsilon}}
\providecommand{\dist}[2]{\ensuremath{\operatorname{d} (#1,#2)}}

\providecommand{\keywords}[1]{
{\let\thefootnote=\relax
\footnote{{\em Keywords}: #1}}
\addtocounter{footnote}{-1}
}

\providecommand{\AMS}[1]{
{\let\thefootnote=\relax
\footnote{{\em AMS Subject Classification} (2010): #1}}
\addtocounter{footnote}{-1}
}

\providecommand{\address}{
{\sc \noindent Department of Mathematics \\
Freie Universit\"at Berlin \\
Arnimallee 6, 14195 Berlin \\
Germany \\}
}

\DeclarePairedDelimiter{\set}{\lbrace}{\rbrace}

\DeclarePairedDelimiter{\abs}{\lvert}{\rvert}
\DeclarePairedDelimiter{\norm}{\lVert}{\rVert}

\theoremstyle{definition}
\newtheorem{definition}{Definition}[section]
\newtheorem*{definition*}{Definition}

\newtheorem*{example*}{Example}

\newtheorem*{remark*}{Remark}
\theoremstyle{plain}
\newtheorem{lemma}[definition]{Lemma}
\newtheorem*{lemma*}{Lemma}
\newtheorem{proposition}[definition]{Proposition}
\newtheorem*{proposition*}{Proposition}

\newtheorem*{theorem*}{Theorem}
\newtheorem{corollary}[definition]{Corollary}
\newtheorem*{corolary*}{Corollary}

\newenvironment{Proof}[1][\proofname]{\begin{proof}[#1] \setlength{\parindent}{0pt}}{\end{proof}}
\newenvironment{Abstract}{\centering\begin{minipage}{0.8\textwidth} \noindent \small {\sc Abstract.}}{\end{minipage}\par}

\usepackage{color}
\definecolor{darkgreen}{rgb}{0,0.5,0}

\numberwithin{equation}{section}
\newtagform{colored}[\color{blue}]{\color{blue}(}{\color{blue})}
\usetagform{colored}

\usepackage[colorlinks,linkcolor=blue,citecolor=red,urlcolor=darkgreen]{hyperref}
\providecommand{\email}{{\it E-mail address:} \href{mailto:hardtke@math.fu-berlin.de}{\tt hardtke@math.fu-berlin.de}}

\usepackage{amsrefs}

\begin{document}

\maketitle

\begin{Abstract}
We study the question whether properties like local/weak almost squareness and local
octahedrality pass down from an absolute sum $X\oplus_F Y$ to the summands $X$ and $Y$.
\end{Abstract}
\AMS{46B20}

\section{Introduction}\label{sec:intro}
First we fix some notation. Throughout this paper we denote by $X$, $Y$, etc. real 
Banach spaces. $X^*$ denotes the dual of $X$, $B_X$ its closed unit ball and $S_X$ its 
unit sphere.\par 
Let us now begin by recalling the following definition (see \cite{godefroy}):
$X$ is called octahedral (OH) if the following holds: for every finite-dimensional subspace 
$F$ of $X$ and every $\eps>0$ there is some $y\in S_X$ such that
\begin{equation*}
\norm{x+y}\geq (1-\eps)(\norm{x}+1) \ \ \forall x\in F.
\end{equation*}
$\ell^1$ is the standard example of an octahedral space. In fact, a Banach space possesses an 
equivalent octahedral norm if and only if it contains an isomorphic copy of $\ell^1$ 
(see \cite{deville}*{Theorem 2.5, p. 106}).\par
In the paper \cite{haller3}, the following weaker forms of octahedrality were introduced: 
$X$ is called locally octahedral (LOH) if for every $x\in X$ and every $\eps>0$ there exists 
$y\in S_X$ such that
\begin{equation*}
\norm{sx+y}\geq (1-\eps)(|s|\norm{x}+1) \ \ \forall s\in \R.
\end{equation*}
$X$ is called weakly octahedral (WOH) if for every finite-dimensional subspace $F$ of $X$, every $x^*\in B_{X^*}$ 
and each $\eps>0$ there is some $y\in S_X$ such that 
\begin{equation*}
\norm{x+y}\geq (1-\eps)(|x^*(x)|+1) \ \ \forall x\in F.
\end{equation*}
The motivation for these definitions was to give dual characterisations of the so called 
diameter-two-properties. For $x^*\in S_{X^*}$ and $\eps>0$, the slice of $B_X$ induced by $x^*$ 
and $\eps$ is the set $S(x^*,\eps):=\set*{z\in B_X,x^*(z)>1-\eps}$. Following the terminology of 
\cite{abrahamsen}, a Banach space $X$ is said to have the local diameter-two-property (LD2P) if every 
slice of $B_X$ has diamter 2 and it is said to have the diameter-two-property (D2P) if every nonempty, 
relatively weakly open subset of $B_X$ has diameter 2. Finally, $X$ is said to have the strong 
diameter-two-property (SD2P) if every convex combination of slices of $B_X$ has diameter 2.\par 
The following results were proved in \cite{haller3}: 
\begin{enumerate}[(a)]
\item $X$ has the LD2P $\iff$ $X^*$ is LOH.
\item $X$ has the D2P $\iff$ $X^*$ is WOH.
\item $X$ has the SD2P $\iff$ $X^*$ is OH.
\end{enumerate}
The result (c) was also proved independently in \cite{becerra-guerrero2}.\par
It is known that the three diameter-two-properties are indeed different. For example, it follows 
from the results on direct sums in \cite{haller3} that $c_0\oplus_2 c_0$ has the D2P but not the SD2P.\par
Concerning the nonequivalence of the LD2P and the D2P, it was proved in \cite{becerra-guerrero4} that every 
Banach space containing an isomorphic copy of $c_0$ can be renormed such that the new space has the LD2P 
but its unit ball contains relatively weakly open subsets of arbitrarily small diameter.\footnote{Note that 
the abbreviation SD2P in \cite{becerra-guerrero4} does not stand for ``strong diameter-two-property'' but 
for ``slice diameter-two-property'', which coincides with the LD2P of \cite{abrahamsen}.}\par 
In \cite{kubiak} it was proved that Ces\`aro function spaces have the D2P.\par 
There are many equivalent formulations of the three octahedrality pro\-perties (see for instance 
\cites{haller2, haller3}). We will recall only those which we need here (they can be found in \cite{haller3}): 
a Banach space $X$ is octahedral if and only if for every $n\in \N$, all $x_1,\dots,x_n\in S_X$ and 
every $\eps>0$ there exists an element $y\in S_X$ such that $\norm{x_i+y}\geq 2-\eps$ for all 
$i=1,\dots,n$. $X$ is locally octahedral if and only if for every $x\in S_X$ and all $\eps>0$
there exists $y\in S_X$ such $\norm{x\pm y}\geq 2-\eps$. We will use these characterisations later without 
further mention.\par
Now we come to the classes of almost square spaces and their relatives. In the paper \cite{abrahamsen2}, 
the following definitions were introduced. A real Banach space $X$ is said to be almost square (ASQ) 
if the following holds: for all $n\in \N$ and all $x_1,\dots,x_n\in S_X$ there exists a sequence 
$(y_k)_{k\in \N}$ in $B_X$ such that $\norm{y_k}\to 1$ and $\norm{x_i+y_k}\to 1$ for all $i=1,\dots,n$.\par 
$X$ is called locally almost square (LASQ) if for every $x\in S_X$ there is a sequence $(y_k)_{k\in \N}$ 
in $B_X$ such that $\norm{y_k}\to 1$ and $\norm{x\pm y_k}\to 1$.\par
$X$ is called weakly almost square (WASQ) if it fulfils the definition of an LASQ space with the 
additional condition that the sequence $(y_k)_{k\in \N}$ converges weakly to zero.\par
Obviously, WASQ implies LASQ. It was shown in \cite{abrahamsen2} that ASQ implies WASQ and that the converse 
of this statement does not hold, while it is not known whether LASQ is strictly weaker than WASQ.\par
The model example of an ASQ space is $c_0$ and it was further proved in \cite{abrahamsen2} that every ASQ space 
contains an isomorphic copy of $c_0$ and, conversely, every separable Banach space containing an isomorphic copy 
of $c_0$ can be equivalently renormed to become ASQ. In \cite{becerra-guerrero3} it was shown that the same holds 
true also for nonseparable spaces. Also, it was proved in \cite{abrahamsen2} that if $X$ is ASQ, than $X^*$ is OH 
(i.\,e. $X$ has the SD2P). By \cite{kubiak}*{Propositions 2.5 and 2.6}, LASQ spaces have the LD2P and WASQ spaces 
have the D2P.\par 
Next we will recall the necessary basics on absolute sums. A norm $F$ on $\R^2$ is called absolute 
if $F(a,b)=F(\abs{a},\abs{b})$ for all $(a,b)\in \R^2$ and it is called normalised if $F(1,0)=1=F(0,1)$. 
If $F$ is an absolute, normalised norm on $\R^2$, and $X$ and $Y$ are two Banach spaces, then the absolute 
sum of $X$ and $Y$ with respect to $F$, denoted by $X\oplus_F Y$, is defined as the direct product $X\times Y$ 
equipped with the norm $\norm{(x,y)}_F=F(\norm{x},\norm{y})$. $X\oplus_F Y$ is again a Banach space.\par 
For every $1\leq p\leq\infty$, the $p$-norm $\norm{\cdot}_p$ on $\R^2$ is an absolute, normalised 
norm and the corresponding sum is just the usual $p$-direct sum of two Banach spaces. We also note 
the following important facts (see for instance \cite{bonsall}*{p. 36, Lemmas 1 and 2}): if $F$ is an 
absolute, normalised norm on $\R^2$, then we have for all $a,b,c,d\in \R$
\begin{enumerate}[1.)]
\item $\abs{a}\leq\abs{c} \ \mathrm{and} \ \abs{b}\leq\abs{d} \ \Rightarrow \ F(a,b)\leq F(c,d)$,
\item $\abs{a}<\abs{c} \ \mathrm{and} \ \abs{b}<\abs{d} \ \Rightarrow \ F(a,b)<F(c,d)$,
\item $\norm{(a,b)}_{\infty}\leq F(a,b)\leq \norm{(a,b)}_1$.
\end{enumerate}
It follows in particular that $\abs{a},\abs{b}\leq F(a,b)$ holds for all $a,b\in \R$.\par
We will also need the following (see \cite{hardtke}): for every $t\in (-1,1)$ there exists a unique 
$f(t)\in (0,1]$ such that $F(t,f(t))=1$. We will call the function $f$ the upper boundary curve
of $B_{(\R^2,F)}$. It is even, concave (hence continuous), decreasing on $[0,1)$ and increasing $(-1,0]$.
Thus it can be extended to a concave, continuous, even function on $[-1,1]$, which will also be denoted by $f$.\par 
Octahedrality properties in $p$-direct sums were already studied in \cite{haller3}. The following
results were proved:
\begin{enumerate}[\upshape(i)]
\item If $X$ or $Y$ is LOH/WOH/OH, then $X\oplus_1 Y$ is LOH/WOH/OH.
\item If $X$ and $Y$ are LOH/WOH, then $X\oplus_pY$ is LOH/WOH for every $1<p\leq\infty$.
\item If $X$ and $Y$ are OH, then $X\oplus_{\infty} Y$ is OH.
\item $X\oplus_p Y$ is never OH for $1<p<\infty$ (provided that $X$ and $Y$ are nontrivial).
\item If $X\oplus_p Y$ is LOH/WOH for $1<p\leq\infty$, then $X$ and $Y$ are LOH/WOH.
\item If $X\oplus_{\infty} Y$ is OH, then $X$ and $Y$ are OH.
\end{enumerate}
Using their duality results, the authors of \cite{haller3} also obtained corresponding results
for diameter-two-properties in $p$-direct sums (see also \cites{abrahamsen, acosta, becerra-guerrero, haller, lopez-perez}
for previous results on diameter-two-properties in $p$-direct sums based on different methods). In \cite{acosta} it was 
also proved that the LD2P and the D2P are stable under all absolute sums, and that $X$ has the LD2P/D2P if $\text{dim}(X)=\infty$,
$X\oplus_F Y$ has the LD2P/D2P for some Banach space $Y$ and $(1,0)$ is an extreme point of $B_{(\R^2,F)}$. Moreover, it was proved in 
\cite{acosta} that $X\oplus_F Y$ does not have the SD2P if $F(1,1)<2$ and $(1,0)$, $(0,1)$ are extreme points of $B_{(\R^2,F)}$ (see 
also \cite{oja}).\par
In the recent paper \cite{haller4} the stability of average roughness (which is a generalisation of octahedrality)
with respect to absolute sums is investigated. The authors of \cite{haller4} also introduce the notion of positive octahedrality 
for an absolute, normalised norm $F$ on $\R^2$, meaning that there exist $c,d\geq 0$ with $F(c,d)=1$ and $F(c+1,d)=F(c,d+1)=2$.
They prove that $X\oplus_F Y$ is octahedral whenever $X$ and $Y$ are octahedral and $F$ is positively octahedral, and, conversely,
if $X\oplus_F Y$ is octahedral for some nontrivial Banach spaces $X$ and $Y$, then $F$ has to be positively octahedral. Analogous 
results for the SD2P in absolute sums are also obtained in \cite{haller4}.\par 
In \cite{abrahamsen2} it is proved that the properties LOH and LASQ are stable under arbitrary (even infinite) 
absolute sums, and that WOH and WASQ are stable under all absolute sums which fulfil a simple density assumption,
including in particular all finite absolute sums.\par 
Also, the following results were obtained in \cite{abrahamsen2} for any two nontrivial Banach spaces $X$ and $Y$.
\begin{enumerate}[\upshape(i)]
\item For $1\leq p<\infty$, $X\oplus_p Y$ is LASQ/WASQ if and only if $X$ and $Y$ are LASQ/WASQ.
\item $X\oplus_{\infty} Y$ is LASQ/WASQ/ASQ if and only if $X$ or $Y$ is LASQ/WASQ/\\ASQ.
\item For $1\leq p <\infty$, $X\oplus_p Y$ is never ASQ.
\end{enumerate}
The purpose of this note is to extend these results by showing that (i) and (iii) also hold if we replace $\norm{\cdot}_p$
by any absolute, normalised norm $F\neq \norm{\cdot}_{\infty}$. We will also prove some results on summands in LOH spaces,
which imply in particular that $X$ and $Y$ are LOH whenever $X\oplus_F Y$ is LOH and $F$ is strictly convex.\par
Finally, we will also discuss some results on ultrapowers of LOH, LASQ, etc. spaces and the closedness of these classes with
respect to the Banach-Mazur distance.

\section{Results and proofs}\label{sec:results}

We start with the following lemma, which is surely well-known, but since the author was not able 
to find it explicitly in the literature, a proof is included here for the reader's convenience.
\begin{lemma}\label{lemma:infty}
Let $F$ be an absolute, normalised norm on $\R^2$.
\begin{enumerate}[\upshape(a)]
\item $F(1,1)=1$ $\Leftrightarrow$ $F=\norm{\cdot}_{\infty}$.
\item $F(1,1)=2$ $\Leftrightarrow$ $F=\norm{\cdot}_1$.
\end{enumerate}
\end{lemma}

\begin{Proof}
(a) Assume that $F(1,1)=1$. Let $(a,b)\in \R^2$ such that $F(a,b)=1$.
Then $\abs{a},\abs{b}\leq 1$. If both $\abs{a}<1$ and $\abs{b}<1$,
then we would have $F(a,b)<F(1,1)=1$ (by the general monotonicity 
properties of absolute norms listed in Section \ref{sec:intro}).
It follows that $\abs{a}=1$ or $\abs{b}=1$, hence $\norm{(a,b)}_{\infty}=1$.\par 
Thus we have $S_{(\R^2,F)}\ssq S_{(\R^2,\norm{\cdot}_{\infty})}$, which implies
$F=\norm{\cdot}_{\infty}$.\par 
(b) Suppose that $F(1,1)=2$, i.\,e. the midpoint of $(0,1)$ and $(1,0)$ lies on 
the unit sphere of $(\R^2,F)$. It follows that the whole line segment from $(0,1)$ 
to $(1,0)$ lies on $S_{(\R^2,F)}$, thus $F(t,1-t)=1$ for every $t\in [0,1]$.\par
Hence we have for every $(a,b)\neq (0,0)$
\begin{equation*}
1=F(\abs{a}/(\abs{a}+\abs{b}),1-\abs{a}/(\abs{a}+\abs{b}))=
F(\abs{a}/(\abs{a}+\abs{b}),\abs{b}/(\abs{a}+\abs{b})),
\end{equation*}
i.\,e. $F(a,b)=\norm{(a,b)}_1$.
\end{Proof}

Before we can come to the first main result on sums of LASQ (etc.) spaces, we have to
prove another auxiliary lemma.

\begin{lemma}\label{lemma:lasq2}
Let $F$ be an absolute, normalised norm on $\R^2$ with $F\neq \norm{\cdot}_{\infty}$
and let $\eps>0$. Then there is a $\delta>0$ such that the following holds:
\begin{equation*}
a,b\geq 0,\ F(a,b)=1 \ \text{and} \ F(a,1)\leq 1+\delta \ \Rightarrow \ b\geq 1-\eps.
\end{equation*}
\end{lemma}

\begin{Proof}
Denote by $f$ the upper boundary curve of $B_{(\R^2,F)}$. If the claim was false, then we could find
sequences $(a_n)_{n\in \N}$, $(b_n)_{n\in \N}$ in $[0,\infty)$ such that $F(a_n,b_n)=1$, 
$F(a_n,1)\leq 1+1/n$ and $b_n<1-\eps$ for each $n\in \N$.\par 
Since $a_n, b_n\leq 1$ for every $n\in \N$, we can find subsequences $(a_{n_k})$, $(b_{n_k})$ such 
that $a_{n_k}\to a$ and $b_{n_k}\to b$ for some $a,b \in [0,1]$. It follows that $F(a,b)=1=F(a,1)$ and 
$b\leq 1-\eps$.\par 
Since $F\neq \norm{\cdot}_{\infty}$ it follows from Lemma \ref{lemma:infty} that $F(1,1)>1$
and hence $a<1$. But then $b=f(a)=1$, be definition of $f$. This is a contradiction since $b<1$.
\end{Proof}

Now we can prove the first main result of this paper.
\begin{proposition}\label{prop:lasq}
If $F$ is any absolute, normalised norm on $\R^2$ with $F\neq \norm{\cdot}_{\infty}$
and $X$ and $Y$ are nontrivial Banach spaces, then the following holds:
\begin{enumerate}[\upshape(i)]
\item If $X\oplus_F Y$ is LASQ, then $X$ and $Y$ are LASQ.
\item If $X\oplus_F Y$ is WASQ, then $X$ and $Y$ are WASQ.
\item $X\oplus_F Y$ is not ASQ.
\end{enumerate}
\end{proposition}
Note that the converses of (i) and (ii) also hold by the general results in \cite{abrahamsen2}.

\begin{Proof}
First we will prove statement (ii). So let $Z:=X\oplus_F Y$ be WASQ and let $y\in S_Y$.
Then there is a weakly null sequence $(z_n=(u_n,v_n))_{n\in \N}$ in $B_Z$ such that 
$\norm{z_n\pm (0,y)}_F\to 1$ and $\norm{z_n}_F\to 1$. Actually, we may assume that 
$\norm{z_n}_F=F(\norm{u_n},\norm{v_n})=1$ for every $n\in \N$.\par 
By Lemma \ref{lemma:lasq2} there exists a sequence $(\delta_n)_{n\in \N}$ in $(0,\infty)$
such that $\delta_n\to 0$ and for every $n\in \N$ the following holds:
\begin{equation*}
a,b\geq 0,\ F(a,b)=1 \ \text{and} \ F(a,1)\leq 1+\delta_n \ \Rightarrow \ b\geq 1-2^{-n}.
\end{equation*}
By passing to a subsequence if necessary, we may assume that $F(\norm{u_n},\norm{y\pm v_n})\leq 1+\delta_n$
for every $n\in \N$. It follows that
\begin{equation*}
F(\norm{u_n},1)\leq \frac{1}{2}(F(\norm{u_n},\norm{y+v_n})+F(\norm{u_n},\norm{y-v_n})\leq 1+\delta_n
\end{equation*}
and hence $\norm{v_n}\geq 1-2^{-n}$ for every $n$.\par 
Since we also have $\norm{v_n}\leq F(\norm{u_n},\norm{v_n})=1$ for each $n$, we obtain $\norm{v_n}\to 1$.
Also, $(v_n)_{n\in \N}$ is a weakly null sequence in $Y$, since $(z_n)_{n\in \N}$ is weakly null in $Z$.\par
We further have $\norm{y\pm v_n}\leq F(\norm{u_n},\norm{y\pm v_n})\leq 1+\delta_n$ and thus
\begin{equation*}
1+\delta_n\geq \norm{y+v_n}\geq 2-\norm{y-v_n}\geq 1-\delta_n \ \ \forall n\in \N,
\end{equation*}
which implies $\norm{y\pm v_n}\to 1$. Thus $Y$ is WASQ.\par 
Since $X\oplus_F Y \cong Y\oplus_{\tilde{F}} X$, where $\tilde{F}(a,b):=F(b,a)$, the same argument also shows 
that $X$ is LASQ. This completes the proof of (ii) and statement (i) is proved analogously.\par
Now we will prove (iii). Assume to the contrary that $X\oplus_F Y$ is ASQ. Since $F\neq \norm{\cdot}_{\infty}$,
we have $F(1,1)>1$ (Lemma \ref{lemma:infty}). Choose $\eps>0$ such that $(1-\eps)F(1,1)>1$.\par 
By Lemma \ref{lemma:lasq2} (applied to $F$ and $\tilde{F}$) there exists a $\delta>0$ such that for all $a,b\geq 0$
with $F(a,b)=1$ the following holds:
\begin{align*}
&F(a,1)\leq 1+\delta \ \Rightarrow b\geq 1-\eps, \\
&F(1,b)\leq 1+\delta \ \Rightarrow a\geq 1-\eps.
\end{align*}
Now let $x\in S_X$ and $y\in S_Y$. Since $X\oplus_F Y$ is ASQ, there exist $u\in X$, $v\in Y$ such that
$F(\norm{u},\norm{v})=1$ and $F(\norm{x\pm u},\norm{v})\leq 1+\delta$, $F(\norm{u},\norm{y\pm v})\leq 1+\delta$.\par
A similar calculation as in the proof of (ii) shows that $F(\norm{u},1)\leq 1+\delta$ and $F(1,\norm{v})\leq 1+\delta$.
It follows that $\norm{u}, \norm{v}\geq 1-\eps$.\par 
But then $1=F(\norm{u},\norm{v})\geq (1-\eps)F(1,1)>1$ and with this contradiction the proof is finished.
\end{Proof}

Next we turn our attention to LOH sums. First recall that a Banach space $X$ is strictly convex (SC) if
$x,y\in S_X$ and $\norm{x+y}=2$ imply $x=y$. The $p$-norms are strictly convex for $1<p<\infty$. We will 
call a point $x\in S_X$ an SC-point of $X$ if $\norm{x+y}<2$ for every $y\in S_X$ with $y\neq x$. Thus $X$ 
is strictly convex if and only if every point of $S_X$ is an SC-point.\par 
Given an absolute, normalised norm $F$ on $\R^2$, set 
\begin{equation*}
r_F:=\inf\set*{a\in [0,1]:\exists\,b\geq 0 \,F(a,b)=1 \ \text{and}\ F(a+1,b)=2}.
\end{equation*}
The following lemma is intuitively clear, but we include a proof for the sake of completeness.

\begin{lemma}\label{lemma:loh1}
Let $F$ be an absolute, normalised norm on $\R^2$ with upper boundary curve $f$. Then
\begin{enumerate}[\upshape(i)]
\item $r_F=1$ \ $\Leftrightarrow$ \ $(1,0)$ is an SC-point of $(\R^2,F)$ or $f(1)>0$;
\item $r_F=0$ \ $\Leftrightarrow$ \ $F=\norm{\cdot}_1$.
\end{enumerate}
\end{lemma}

\begin{Proof}
(i) If $r_F<1$, then there must be $a\in [0,1)$ and $b>0$ such that $F(a,b)=1$
and $F(a+1,b)=2$. Hence $(1,0)$ is not an SC-point of $(\R^2,F)$.\par 
Moreover, the whole line segment from $(1,0)$ to $(a,b)$ belongs to the unit
sphere of $(\R^2,F)$, which implies that $f(a+t(1-a))=(1-t)b$ for $t\in [0,1)$.
Hence $f(s)=\frac{1-s}{1-a}b$ for $s\in [a,1)$. Thus $f(1)=0$. This shows ``$\Leftarrow$''.\par
Now assume that $r_F=1$ and $(1,0)$ is not an SC-point of $(\R^2,F)$. Then we can find 
$(a,b)\in \R^2$ such that $F(a,b)=1$ and $F(1+a,b)=2$ but $(a,b)\neq (1,0)$. Without loss
of generality we may assume that $a,b>0$. Since $r_F=1$ it follows that $a=1$, hence $F(1,b)=1$.\par
If $f(1)=0$ there would be $s\in [0,1)$ such that $f(s)<b$. But then we obtain a contradiction since 
$1=F(s,f(s))<F(1,b)=1$. Thus we must have $f(1)>0$.\par 
(ii) Suppose that $r_F=0$. This easily implies $F(1,1)=2$ and thus we have
$F=\norm{\cdot}_1$ by Lemma \ref{lemma:infty}. The converse is clear.
\end{Proof}

We need two more auxiliary lemmas.
\begin{lemma}\label{lemma:loh2}
Let $F$ be an absolute, normalised norm on $\R^2$. If $a,b\geq 0$ such that $F(a,b)=1$ 
and $0\leq c\leq 1+a$ such that $F(c,b)=2$, then $c=1+a$ and $a\geq r_F$.
\end{lemma}

\begin{Proof}
First note that under the above assumptions we have $2=F(c,b)\leq F(1+a,b)\leq F(1,0)+F(a,b)=2$. 
Hence, by definition of $r_F$, we must have $a\geq r_F$.\par  
Now we denote again by $f$ the upper boundary curve of $F$ and distinguish two cases.\par
Case 1: $r_F=1$. Then $a=1$ and $1<c\leq 2$ (if $c\leq 1$ we would obtain $2=F(c,b)\leq F(1,b)=1$).
Thus $a^{\prime}:=c-1\in (0,1]$ and $F(1+a^{\prime},b)=2$ as well as 
$1=F(1,b)\geq F(a^{\prime},b)\geq F(c,b)-F(1,0)=1$.\par 
Since $r_F=1$ it follows that $a^{\prime}=1$, i.\,e. $c=2$.\par 
Case 2: $r_F<1$. Then we have $r_F\leq a<c\leq 1+a$. Put $d:=c-a\in (0,1]$ and let $w>0$ such that 
$F(r_F,w)=1$ and $F(1+r_F,w)=2$.\par 
Then the line segment from $(1,0)$ to $(r_F,w)$ lies completely in $S_{(\R^2,F)}$ and we obtain
\begin{equation*}
f(s)=h(s-1) \ \ \forall s\in [r_F,1], \ \text{where}\ h:=\frac{w}{r_F-1}.
\end {equation*}
We also put $g(s):=(d+1)f(s/(d+1))$ for $s\in [0,d+1]$. Then $F(s,g(s))=d+1$ for every such $s$.\par 
It is easy to see that $c/(d+1)\in [a,1]\ssq [r_F,1]$ and thus we have $g(c)=(d+1)f(c/(d+1))=h(c-(d+1))
=h(a-1)$.\par 
If $a<1$ we have $f(a)=b$ and since $a\geq r_F$ it follows that $f(a)=h(a-1)$, thus $g(c)=b$.\par 
For $a=1$ we must have $b=0=g(c)$ (otherwise there is $s\in [0,1)$ such that $f(s)<b$ and hence
$1=F(s,f(s))<F(1,b)=F(a,b)=1$).\par 
Thus we always have $g(c)=b$, which imples that $2=F(c,b)=F(c,g(c))=d+1$. Hence $d=1$ and $c=1+a$.
\end{Proof}

\begin{lemma}\label{lemma:loh3}
Let $F$ be an absolute, normalised norm on $\R^2$ and let $\eps>0$. Then there exists $\delta>0$ 
such that the following holds: whenever $a,b\geq 0$ with $F(a,b)=1$ and $0\leq c\leq 1+a$ with 
$F(c,b)\geq 2-\delta$, then $c\geq 1+r_F-\eps$.
\end{lemma}

\begin{Proof}
This follows from Lemma \ref{lemma:loh2} and a standard compactness argument.
\end{Proof}

Now we can prove the main result on local octahedrality in absolute sums.
\begin{proposition}\label{prop:loh}
Let $F$ be an absolute, normalised norm on $\R^2$ and $X$ and $Y$ nontrivial Banach spaces such that 
$X\oplus_F Y$ is LOH. Then for every $x\in S_X$ and every $\eps>0$ there is a $z\in S_X$ such 
that $\norm{x\pm z}\geq 2r_F-\eps$.
\end{proposition}

\begin{Proof}
Let $x\in S_X$ and $\eps>0$. Choose $\delta>0$ according to Lemma \ref{lemma:loh3}
for the parameter $\eps/2$. Since $X\oplus_F Y$ is LOH, we can find $(u,v)\in S_{X\oplus_F Y}$ 
such that $F(\norm{x\pm u},\norm{v})\geq 2-\delta$.\par
Because of $\norm{x\pm u}\leq 1+\norm{u}$ and $F(\norm{u},\norm{v})=1$ this implies
$\norm{x\pm u}\geq 1+r_F-\eps/2$.\par 
It follows that $\norm{u}\geq \norm{x+u}-1\geq r_F-\eps/2$.\par 
Now put $z:=u/\norm{u}$. Then 
\begin{align*}
\norm{x\pm z}\geq \norm{x\pm u}-\norm{u-u/\norm{u}}\geq 1+r_F-\frac{\eps}{2}-(1-\norm{u})\geq 2r_F-\eps.
\end{align*}
\end{Proof}

It follows in particular from Proposition \ref{prop:loh} that $X$ is LOH if $X\oplus_F Y$ is LOH and $r_F=1$
(which, by Lemma \ref{lemma:loh1}, is equivalent to the fact that $f(1)>0$ or $(1,0)$ is an SC-point of $(\R^2,F)$).\par 
More generally, for any Banach space $X$ we may define 
\begin{equation*}
s(X):=\sup\set*{s\in [0,2]:\forall x\in S_X\,\forall \eps>0\,\exists y\in S_X \,\norm{x\pm y}\geq s-\eps}.
\end{equation*}
Then $X$ is LOH if and only if $s(X)=2$ and Proposition \ref{prop:loh} reads: if $X\oplus_F Y$ is LOH, 
then $s(X)\geq 2r_F$.\par 
Note that $s(\R)=0$, while it easily follows from Riesz's Lemma that $s(X)\geq 1$ whenever $\text{dim}(X)\geq 2$.
The following statements are also easy to verify: $s(\ell^{\infty})=s(c_0)=1=s(\ell^{\infty}(n))$ for $n\geq 2$,
where $\ell^{\infty}$ is the space of bounded sequences, $c_0$ is the space of null sequences (both equipped with 
the supremum norm), and $\ell^{\infty}(n)$ is the space $\R^n$ equipped with the maximum norm. It is also easy to 
prove that $s(H)=\sqrt{2}$ for any Hilbert space $H$ with $\text{dim}(H)\geq 2$.\par 
Putting everything together, we obtain the following corollary (for (b) note that $r_F>0$ if $F\neq \norm{\cdot}_1$
(Lemma \ref{lemma:loh1})).

\begin{corollary}\label{cor:loh}
Let $F$ be an absolute, normalised norm on $\R^2$ and $X$ and $Y$ nontrivial Banach spaces.
Then the following holds:
\begin{enumerate}[\upshape(a)]
\item If $r_F=1$ and $X\oplus_F Y$ is LOH, then so is $X$. In particular, this holds if $F$ is strictly convex
or $F=\norm{\cdot}_{\infty}$.
\item If $F\neq \norm{\cdot}_1$, then $\R\oplus_F Y$ is not LOH.
\item If $r_F>1/2$, then $\ell^{\infty}\oplus_F Y$, $c_0\oplus_F Y$ and $\ell^{\infty}(n)\oplus_F Y$ for $n\in \N$
are not LOH.
\item If $r_F>1/\sqrt{2}$ and $H$ is a Hilbert space, then $H\oplus_F Y$ is not LOH.
\end{enumerate}
\end{corollary}

Of course, it is also possible to prove results analogous to Proposition \ref{prop:loh} and Corollary \ref{cor:loh}
for the second summand by modifying the definition of $r_F$ accordingly (i.\,e. using instead $r_{\tilde{F}}$, where
$\tilde{F}(a,b):=F(b,a)$).\par 
The author does not know whether there are any analogous results for WOH spaces, but let us remark that the above 
proof-techniques could also be used to show that $X\oplus_F Y$ is not octahedral if $X$ and $Y$ are nontrivial Banach 
spaces, $F\neq \norm{\cdot}_{\infty}$ and $r_F=1=r_{\tilde{F}}$. However, this result already follows from the more general 
results on octahedrality in absolute sums that were proved in \cite{haller4} and that we have already mentioned in the introduction.\par

Let us note one more corollary concerning the LD2P (recall that a norm is smooth if it is G\^ateaux-differentiable at each nonzero point).
\begin{corollary}\label{cor:ld2p}
If $F$ is a smooth, absolute, normalised norm on $\R^2$ and $X$ and $Y$ are nontrivial Banach spaces such that $X\oplus_F Y$
has the LD2P, then $X$ has the LD2P.
\end{corollary}

\begin{Proof}
It is well-known that a finite-dimensional Banach space is smooth if and only if its dual is strictly convex.
If we put 
\begin{equation*}
F^*(c,d):=\sup\set*{|ac|+|bd|:(a,b)\in B_{(\R^2,F)}},
\end{equation*}
then $F^*$ is an absolute, normalised norm on $\R^2$ and $(X\oplus_F Y)^*\cong X^*\oplus_{F^*} Y^*$ (this is a standard 
fact from the theory of absolute sums, which is easy to prove). The claim now follows from Corollary \ref{cor:loh} and 
the duality between LOH and LD2P (\cite{haller3}).
\end{Proof}

Next we will consider ultrapowers of OH/LOH and ASQ/LASQ spaces. First we recall the necessary definitions.
Given a free ultrafilter $\U$ on $\N$ and a bounded sequence $(a_n)_{n\in \N}$ of real numbers, there exists 
(by a compactness argument) a unique number $a\in \R$ such that for every $\eps>0$ we have $\set*{n\in \N:\abs{a_n-a}<\eps}\in\U$. 
It is called the limit of $(a_n)_{n\in \N}$ along $\U$ and will be denoted by $\lim_{n, \U}a_n$.\par
For a Banach space $X$, denote by $\ell^{\infty}(X)$ the space of all bounded sequences in $X$ and set 
$\mathcal{N}_{\U}:=\set*{(x_n)\in \ell^{\infty}(X):\lim_{n,\U}\norm{x_n}=0}$. The ultarpower $X^{\U}$ of $X$ 
with respect to $\U$ is the quotient space $\ell^{\infty}(X)/\mathcal{N}_{\U}$ equipped with the (well-defined)
norm $\norm{[(x_n)]}_{\U}:=\lim_{n,\U}\norm{x_n}$. $X^{\U}$ is again a Banach space (for more information 
on ultraproducts see for example \cite{heinrich}).\par 
We have the following observations concerning octahedrality in ultrapowers.

\begin{proposition}\label{prop:ultraocta}
Let $X$ be a Banach space and $\U$ a free ultrafilter on $\N$. Then the following assertions
are equivalent:
\begin{enumerate}[\upshape(i)]
\item $X$ is octahedral.
\item For all $z_1,\dots,z_n\in S_{X^{\U}}$ there exists an element $z\in S_{X^{\U}}$ such
that $\norm{z_i+z}_{\U}=2$ for every $i\in \set*{1,\dots,n}$.
\item $X^{\U}$ is octahedral.
\end{enumerate}
Likewise, the following statements are equivalent:
\begin{enumerate}[\upshape(i)]
\item $X$ is locally octahedral.
\item For every $z\in S_{X^{\U}}$ there is some $\tilde{z}\in S_{X^{\U}}$ such that 
$\norm{z\pm \tilde{z}}_{\U}=2$.
\item $X^{\U}$ is locally octahedral.
\end{enumerate}
\end{proposition}

\begin{Proof}
We will only prove the statement for octahedral spaces. The proof for local octahedrality is completely analogous.\par 
So let us first assume that $X$ is OH and let $z_1,\dots,z_n\in S_{X^{\U}}$. Let $(x_{i,k})_{k\in \N}$ be a representative
of $z_i$. We may assume that $x_{i,k}\neq 0$ for all $k\in \N$ and all $i=1,\dots,n$.\par
Since $X$ is octahedral we can find, for each $k\in \N$, an element $x_k\in S_X$ such that
\begin{equation*}
\norm*{\frac{x_{i,k}}{\norm{x_{i,k}}}+x_k}\geq 2-2^{-k} \ \ \forall i=1,\dots,n.
\end{equation*}
Let $z:=[(x_k)_{k\in \N}]\in S_{X^{\U}}$. For each $k\in \N$ and every $i\in \{1,\dots,n\}$ we have
\begin{equation*}
\norm{x_{i,k}+x_k}\geq\norm*{\frac{x_{i,k}}{\norm{x_{i,k}}}+x_k}-\norm*{\frac{x_{i,k}}{\norm{x_{i,k}}}-x_{i,k}}\geq 2-2^{-k}-|1-\norm{x_{i,k}}|.
\end{equation*}
Since $\lim_{k,\U}\norm{x_{i,k}}=\norm{z_i}_{\U}=1$ it follows that $\norm{z_i+z}_{\U}=\lim_{k,\U}\norm{x_{i,k}+x_k}=2$
for all $i\in \{1,\dots,n\}$. This proves (i) $\Rightarrow$ (ii).\par
(ii) $\Rightarrow$ (iii) is clear.\par
(iii) $\Rightarrow$ (i): Let $x_1,\dots,x_n\in S_X$ and $\eps>0$. We consider $X$ as a subspace of $X^{\U}$ (via the canonical embedding).
Since $X^{\U}$ is octahedral, there exists $y=[(y_k)_{k\in \N}]\in S_{X^{\U}}$ such that 
\begin{equation*}
\lim_{k,\U}\norm{x_i+y_k}=\norm{x_i+y}_{\U}\geq 2-\eps \ \ \forall i=1,\dots,n.
\end{equation*}
It follows that $B_i:=\{k\in\N:\norm{x_i+y_k}\geq 2-2\eps\}\in\U$ for all $i=1,\dots,n$. Since $\lim_{k,\U}\norm{y_k}=\norm{y}_{\U}=1$
we also have $A:=\{k\in \N:|\norm{y_k}-1|\leq\eps\}\in \U$. Hence $M:=A\cap B_1\cap\dots\cap B_n\in \U$ and in particular, $M\neq\emptyset$.\par
Now let $k_0\in M$. Then we have, for each $i\in \{1,\dots,n\}$,
\begin{equation*}
\norm*{x_i+\frac{y_{k_0}}{\norm{y_{k_0}}}}\geq\norm{x_i+y_{k_0}}-\norm*{y_{k_0}-\frac{y_{k_0}}{\norm{y_{k_0}}}}\geq2-2\eps-|\norm{y_{k_0}}-1|\geq2-3\eps,
\end{equation*}
thus $X$ is octahedral.
\end{Proof}

For weakly octahedral spaces, the situation seems to be more complicated. Let us first introduce one more notation:
if $s=(x_n^*)_{n\in \N}$ is a sequence in $S_{X^*}$, then we may define a norm-one functional $\varphi_s$ 
on $X^{\U}$ by $\varphi_s([(x_n)]):=\lim_{n,\U}x_n^*(x_n)$.\par 
Using the characterisation for WOH spaces from \cite{haller3}*{Proposition 2.2}, one can easily prove the following:
if $X$ is WOH, then for all $z_1,\dots,z_n\in S_{X^{\U}}$ and every sequence $s=(x_n^*)_{n\in \N}$ in $S_{X^*}$
there exists $z\in S_{X^{\U}}$ such that
\begin{equation*}
\norm{z_i+tz}_{\U}\geq \abs{\varphi_s(z_i)}+t \ \ \forall i\in \set*{1,\dots,n},\ \forall t>0.
\end{equation*}
However, it is not clear whether the converse of this statement also holds, nor whether it is equivalent to 
the weak octahedrality of $X^{\U}$.\par
Similarly to Proposi\-tion \ref{prop:ultraocta} one can also prove the following result for ASQ/LASQ spaces 
(we skip the details).
\begin{proposition}\label{prop:ultraasq}
Let $X$ be a Banach space and $\U$ a free ultrafilter on $\N$. Then the following assertions
are equivalent:
\begin{enumerate}[\upshape(i)]
\item $X$ is ASQ.
\item For all $z_1,\dots,z_n\in S_{X^{\U}}$ there exists an element $z\in S_{X^{\U}}$ such
that $\norm{z_i+z}_{\U}=1$ for every $i\in \set*{1,\dots,n}$.
\item $X^{\U}$ is ASQ.
\end{enumerate}
Likewise, the following statements are equivalent:
\begin{enumerate}[\upshape(i)]
\item $X$ is LASQ.
\item For every $z\in S_{X^{\U}}$ there is some $\tilde{z}\in S_{X^{\U}}$ such that 
$\norm{z\pm \tilde{z}}_{\U}=1$.
\item $X^{\U}$ is LASQ.
\end{enumerate}
\end{proposition}

For WASQ spaces, the situation is again a bit more involved. First we note the following equivalent characterisation 
of WASQ spaces with separable dual (the proof is easy and will therefore be omitted).
\begin{lemma}\label{lemma:WASQ}
Let $X$ be a Banach space. If $X$ is WASQ, then the following holds: for every $x\in S_X$, every $\eps>0$ and all
$x_1^*,\dots,x_n^*\in S_{X^*}$ there exists a $y\in S_X$ such that $\norm{x\pm y}\leq 1+\eps$ and $x_i^*(y)\leq \eps$
for every $i=1,\dots,n$.\par 
If $X^*$ is separable, then the converse of this statement also holds.
\end{lemma}

Using this lemma, it is easy to show the next result (again the details are omitted).
\begin{proposition}\label{prop:ultraWASQ}
Let $X$ be a Banach space and $\U$ a free ultrafilter on $\N$. If $X$ is WASQ, then for every $z\in S_{X^{\U}}$
and all double-sequences $(x_{ik}^*)_{i,k\in \N}$ in $S_{X^*}$ there is some $\tilde{z}\in S_{X^{\U}}$ satisfying
$\norm{z\pm \tilde{z}}_{\U}=1$ and $\varphi_{s_i}(\tilde{z})=0$ for every $i\in \N$, where $s_i:=(x_{ik}^*)_{k\in \N}$.\par
If $X^*$ is separable, then the converse also holds.
\end{proposition}

Again, it is not clear whether the property in Proposition \ref{prop:ultraWASQ} is equivalent to the weak almost 
squareness of $X^{\U}$.\par 
Finally, we would like to show that the classes of OH/WOH/LOH/ASQ/\\LASQ spaces and the classes of spaces 
with the SD2P/D2P/LD2P are closed with respect to the Banach-Mazur distance. Recall that this distance
between two isomorphic Banach spaces $X$ and $Y$ is defined by
\begin{equation*}
\dist{X}{Y}:=\inf\set*{\norm{T}\norm{T^{-1}}:T:X \rightarrow Y \ \text{is an isomorphism}}.
\end{equation*}

\begin{proposition}\label{prop:BanachMazur}
Let $X$ be a Banach space such that for every $\delta>0$ there is some OH/WOH/LOH/ASQ/LASQ space 
$Y$ isomorphic to $X$ with $\dist{X}{Y}<1+\delta$. Then $X$ is also OH/WOH/LOH/ASQ/LASQ.
\end{proposition}

\begin{Proof}
The proofs are all similar, so we will only show the most complicated case of WOH spaces.
Let $x_1,\dots,x_n\in S_X$, $x^*\in B_{X^*}$ and $0<\eps<1/2$. Choose $0<\delta<\eps^2$
such that $\sqrt{\delta}/(1-2\sqrt{\delta})\leq \eps$.\par 
By assumption there is a WOH space $Y$ isomorphic to $X$ such that $\dist{X}{Y}<1+\delta$.
Hence we can find an isomorphism $T:X \rightarrow Y$ such that $\norm{T}=1$ and $\norm{T^{-1}}<1+\delta$.\par
Put $y_i:=Tx_i\in B_Y\sm \set*{0}$ for $i=1,\dots,n$ and $y^*:=(T^*)^{-1}x^*=x^*\circ T^{-1}\in (1+\delta)B_{Y^*}$.\par
Since $Y$ is WOH, there exists, by \cite{haller3}*{Proposition 2.2}, an element $y\in S_Y$ such that
\begin{equation}\label{eq:BanachMazur1}
\norm{y_i+ty}\geq (1-\delta)(\abs{y^*(y_i)}/(1+\delta)+t) \ \ \forall i\in \set*{1,\dots,n},\ \forall t>0.
\end{equation}
Let $z:=T^{-1}y/\norm{T^{-1}y}\in S_X$. Then we have for each $i\in \set*{1,\dots,n}$ and every $t>0$
\begin{equation*}
\norm{x_i+tz}\geq \norm{x_i+tT^{-1}y}-t\abs{\norm{T^{-1}y}-1}\geq \norm{y_i+ty}-t\delta.
\end{equation*}
Combining this with \eqref{eq:BanachMazur1} and observing that $y^*(y_i)=x^*(x_i)$ we obtain
\begin{equation*}
\norm{x_i+tz}\geq \norm{(1+\delta)y_i+ty}-\delta-t\delta\geq (1-\delta)(\abs{x^*(x_i)}+t)-(1+t)\delta
\end{equation*}
for every $i\in \set*{1,\dots,n}$ and every $t>0$.\par 
Now if $t\geq \eps$, then by the choice of $\delta $ we obtain $(\sqrt{\delta}-2\delta)t\geq \delta$, which 
implies $(1-2\delta)t-\delta\geq (1-\sqrt{\delta})t$.\par
Thus for every $t\geq \eps$ and every $i\in \set*{1,\dots,n}$ we have
\begin{equation*}
\norm{x_i+tz}\geq (1-\delta)\abs{x^*(x_i)}+(1-\sqrt{\delta})t\geq (1-\eps)(\abs{x^*(x_i)}+t).
\end{equation*}
By \cite{haller3}*{Proposition 2.2} this implies that $X$ is WOH.
\end{Proof}

Proposition \ref{prop:BanachMazur} and the duality between diameter-two- and octahedrality properties,
together with the fact that $\dist{X^*}{Y^*}\leq \dist{X}{Y}$ holds for all Banach spaces $X$ and $Y$,
immediately yield the following corollary.
\begin{corollary}\label{cor:BanachMazur}
Let $X$ be a Banach space such that for every $\delta>0$ there is a Banach space $Y$ isomorphic to $X$ 
which has the SD2P/D2P/LD2P and satisfies $\dist{X}{Y}<1+\delta$. Then $X$ also has the SD2P/D2P/LD2P.
\end{corollary}

Concerning WASQ spaces, using Lemma \ref{lemma:WASQ} one can show the following result (again an easy 
proof is omitted).
\begin{proposition}\label{prop:BanachMazurWASQ}
Let $X$ be a Banach space such that for every $\delta>0$ there is some WASQ space $Y$ isomorphic to 
$X$ with $\dist{X}{Y}<1+\delta$. If $X^*$ is separable, then $X$ is also WASQ.
\end{proposition}

\address
\email


\begin{bibdiv}
\begin{biblist}

\bib{abrahamsen}{article}{
  title={Remarks on diameter 2 properties},
  author={Abrahamsen, T.},
  author={Lima, V.},
  author={Nygaard, O.},
  journal={J. Conv. Anal.},
  volume={20},
  date={2013},
  pages={439--452}
  }

\bib{abrahamsen2}{article}{
  title={Almost square Banach spaces},
  author={Abrahamsen, T. A.},
  author={Langemets, J.},
  author={Lima, V.},
  journal={J. Math. Anal. Appl.},
  volume={434},
  number={2},
  date={2016},
  pages={1549--1565}
  }
  
\bib{acosta}{article}{
  title={Stability results of diameter two properties},
  author={Acosta, M. D.},
  author={Becerra Guerrero, J.},
  author={L\'opez-P\'erez, G.},
  journal={J. Conv. Anal.},
  volume={22},
  number={1},
  date={2015},
  pages={1--17}
  }
  
\bib{becerra-guerrero}{article}{
  title={Relatively weakly open subsets of the unit ball in function spaces},
  author={Becerra Guerrero, J.},
  author={L\'opez-P\'erez, G.},
  journal={J. Math. Anal. Appl.},
  volume={315},
  date={2006},
  pages={544--554}
  }

\bib{becerra-guerrero2}{article}{
  title={Octahedral norms and convex combination of slices in Banach spaces},
  author={Becerra Guerrero, J.},
  author={L\'opez-P\'erez, G.},
  author={Rueda Zoca, A.},
  journal={J. Funct. Anal.},
  volume={266},
  number={4},
  date={2014},
  pages={2424--2435}
  }
  
\bib{becerra-guerrero4}{article}{
  title={Big slices versus big relatively weakly open subsets in Banach spaces},
  author={Becerra Guerrero, J.},
  author={L\'opez-P\'erez, G.},
  author={Rueda Zoca, A.},
  journal={J. Math. Anal. Appl.},
  volume={428},
  date={2015},
  pages={855--865}
  }
  
\bib{becerra-guerrero3}{article}{
  title={Some results on almost square Banach spaces},
  author={Becerra Guerrero, J.},
  author={L\'opez-P\'erez, G.},
  author={Rueda Zoca, A.},
  journal={J. Math. Anal. Appl.},
  volume={438},
  number={2},
  date={2016},
  pages={1030--1040}
  }
    
\bib{bonsall}{book}{
  title={Numerical Ranges II},
  author={Bonsall, F. F.},
  author={Duncan, J.},
  series={London Math. Soc. Lecture Note Series},
  volume={10},
  publisher={Cambridge University Press},
  address={Cambridge},
  date={1973}
  }
  
\bib{deville}{book}{
  title={Smoothness and renormings in Banach spaces},
  author={Deville, R.},
  author={Godefroy, G.},
  author={Zizler, V.},
  publisher={Longman Scientific \& Technical},
  address={Harlow},
  series={Pitman Monographs and Surveys in Pure and Applied Mathematics},
  volume={64},
  date={1993}
  }
  
\bib{godefroy}{article}{
  title={Metric characterization of first Baire class linear forms and octahedral norms},
  author={Godefroy, G.},
  journal={Studia Math.},
  volume={95},
  number={1},
  date={1989},
  pages={1--15}
  }
  
\bib{haller}{article}{
  title={Two remarks on diameter 2 properties},
  author={Haller, R.},
  author={Langemets, J.},
  journal={Proc. Estonian Acad. Sci.},
  volume={63},
  number={1},
  date={2014},
  pages={2--7}
  }
  
\bib{haller2}{article}{
  title={Geometry of Banach spaces with an octahedral norm},
  author={Haller, R.},
  author={Langemets, J.},
  journal={Acta Comment. Univ. Tartuensis Math.},
  volume={18},
  number={1},
  date={2014},
  pages={125-133}
  }
  
\bib{haller3}{article}{
  title={On duality of diameter 2 properties},
  author={Haller, R.},
  author={Langemets, J.},
  author={P\~{o}ldvere, M.},
  journal={J. Conv. Anal.},
  volume={22},
  number={2},
  date={2015},
  pages={465--483}
  }
  
\bib{haller4}{article}{
  title={Stability of average roughness, octahedrality, 
  and strong diameter two properties of Banach spaces with respect to absolute sums},
  author={Haller, R.},
  author={Langemets, J.},
  author={Nadel, R.},
  date={2017},
  pages={19p.},
  note={Preprint, available at \url{http://arxiv.org/abs/1702.03140}}
  }
  
\bib{hardtke}{article}{
  title={Ball generated property of direct sums of Banach spaces}, 
  author={Hardtke, J.-D.},
  journal={Math. Proc. Royal Irish Acad.},
  volume={115A},
  number={2},
  date={2015},
  pages={137--144}
  }
   
\bib{heinrich}{article}{
  title={Ultraproducts in Banach space theory},
  author={Heinrich, S.},
  journal={J. Reine Angew. Math.}, 
  volume={313},
  date={1980},
  pages={72--104}
  }
  
\bib{kubiak}{article}{
  title={Some geometric properties of the Ces\`aro function spaces},
  author={Kubiak, D.},
  journal={J. Convex Anal.},
  volume={21},
  number={1},
  date={2014},
  pages={189--200}
}
  
\bib{lopez-perez}{article}{
  title={The big slice phenomena in $M$-embedded and $L$-embedded spaces},
  author={L\'opez-P\'erez, G.},
  journal={Proc. Amer. Math. Soc.},
  volume={134},
  date={2005},
  pages={273--282}
  }
  
\bib{oja}{article}{
  title={Sums of slices in direct sums of Banach spaces},
  author={Oja, E.},
  journal={Proc. Estonian Acad. Sci.},
  volume={63},
  number={1},
  date={2014},
  pages={8--10}
  }

\end{biblist}
\end{bibdiv}
\end{document}